\documentclass[11pt,a4paper]{amsart}

\usepackage{amsmath,amssymb}
\newtheorem{theorem}{Theorem}[section]

\newtheorem{proposition}{Proposition}[section]
\newtheorem{definition}{Definition}[section]
\newtheorem{remark}{Remark}[section]
\numberwithin{equation}{section}
\newtheorem{example}{Example}[section]

\begin{document}

\title[Nonuniform Stability for Skew-Evolution Semiflows]
{On nonuniform exponential stability for skew-evolution semiflows
on Banach spaces}
\author{Codru\c{t}a Stoica and Mihail Megan}

\date{}

\begin{abstract}
The paper considers some concepts of nonuniform asymptotic
stability for skew-evolution semiflows on Banach spaces. The
obtained results clarify differences between the uniform and
nonuniform cases. Some examples are included to illustrate the
results.
\end{abstract}



\maketitle

\section{Introduction}\label{Intro}

The exponential stability plays a central role in the theory of
asymptotic behaviors for dynamical systems. In this paper we
consider the more general concepts of nonuniform exponential
stability for skew-evolution semiflows on Banach spaces. These
seem to be more appropriate for the study of evolution equations
in the nonuniform case, because of the fact that they depend on
three variables, contrary to a skew-product semiflow or an
evolution operator, which depend only on two, and, hence, the
study of asymptotic behaviors for skew-evolution semiflows in the
nonuniform setting arises as natural, relative to the third
variable.

Our main objectives are to establish relations between these
concepts and to give some integral characterizations for them. We
also remark that we use the concept of nonuniform exponential
stability, given and studied in the papers of L. Barreira and C.
Valls, as for example \cite{BaVa_LNM}, \cite{BaVa_NA09} or
\cite{BaVa_NA10}, and which we call \textit{"Barreira-Valls
exponential stability"}.

The paper presents some generalizations for the results obtained
in the uniform case in our paper \cite{StMe_NA}.

We remark that Theorems \ref{Datko} and \ref{Datko_BV} are
generalizations of Datko type for the nonuniform exponential
stability in the case of skew-evolution semiflows. The uniform
exponential stability was characterized by R. Datko in
\cite{Da_JMA}. The particular case of evolution operators was
considered by C. Bu\c{s}e in \cite{Bu_RSMUPT} and by M. Megan,
A.L. Sasu and B. Sasu in \cite{MeSaSa_MIR}. Theorem \ref{Rolewicz}
is the nonuniform variant for skew-evolution semiflows of the
known result of S. Rolewicz in \cite{Ro_JMAA}. Theorem
\ref{Barbashin} is a generalization of a result proved by E.A.
Barbashin in \cite{Ba_NAU}. A similar result was obtained
Bu\c{s}e, M. Megan, M. Prajea and P. Preda for the uniform
exponential stability in \cite{BuMePrPr_IEOT}.

Some illustrating examples clarify the connections between the
stability concepts considered in this paper.

\section{Skew-evolution semiflows}\label{Def_ses}

Let $( X,d)$ be a metric space, $ V$ a Banach space and $ V^{*}$
its topological dual. Let $\mathcal{B}( V)$ be the space of all
$V$-valued bounded operators defined on $ V$. The norm of vectors
on $ V$ and on $ V^{*}$ and of operators on $\mathcal{B}( V)$ is
denoted by $\left\Vert \cdot \right\Vert$. Let us consider $Y=
X\times V$ and $T=\left\{ (t,t_{0})\in \mathbb{R}_{+}^{2}:t\geq
t_{0}\right\}$. $I$ is the identity operator.

\begin{definition}
A mapping $\varphi: T\times  X\rightarrow  X$ is called
\emph{evolution semiflow} on $ X$ if the following
propositionerties are satisfied:

$(es_{1})$ $\varphi(t,t,x)=x, \ \forall (t,x)\in
\mathbb{R}_{+}\times  X$;

$(es_{2})$ $\varphi(t,s,\varphi(s,t_{0},x))=\varphi(t,t_{0},x), \
\forall (t,s),(s,t_{0})\in  T, \ \forall x\in
 X$.
\end{definition}

\begin{definition}
A mapping $\Phi: T\times  X\rightarrow \mathcal{B}( V)$ is called
\emph{evolution cocycle} over an evolution semiflow $\varphi$ if
it satisfies following propositionerties:

$(ec_{1})$ $\Phi(t,t,x)=I, \ \forall t\geq0,\ \forall x\in
 X$;

$(ec_{2})$
$\Phi(t,s,\varphi(s,t_{0},x))\Phi(s,t_{0},x)=\Phi(t,t_{0},x),
\forall (t,s), (s,t_{0})\in  T,\forall x\in  X$.
\end{definition}

If $\Phi$ is an evolution cocycle over an evolution semiflow
$\varphi$, then the mapping
\begin{equation}
C: T\times Y\rightarrow Y, \
C(t,s,x,v)=(\varphi(t,s,x),\Phi(t,s,x)v)
\end{equation}
is called \emph{skew-evolution semiflow} on $ Y$.

\begin{remark}\rm
The concept of skew-evolution semiflow generalizes the notion of
skew-product semiflow, considered and studied by M. Megan, A.L.
Sasu and B. Sasu in \cite{MeSaSa_BBMS} and \cite{MeSaSa_MB}, where
the mappings $\varphi$ and $\Phi$ do not depend on the variables
$t$ and $x$.
\end{remark}

\begin{example}\rm
Let $ X=\mathbb{R}_{+}$. The mapping $\varphi:
T\times\mathbb{R}_{+} \rightarrow \mathbb{R}_{+}, \
\varphi(t,s,x)=t-s+x$ is an evolution semiflow on
$\mathbb{R}_{+}$. For every evolution operator $E: T\rightarrow
\mathcal{B}(V)$ (i.e. $E(t,t)=I$, $ \forall t\in \mathbb{R}_{+}$
and $E(t,s)E(s,t_{0})=E(t,t_{0})$, $\forall (t,s),(s,t_{0})\in T$)
we obtain that $\Phi_{E}: T\times \mathbb{R}_{+}\rightarrow
\mathcal{B}(V), \ \Phi_{E}(t,s,x)=E(t-s+x,x)$ \noindent is an
evolution cocycle on $V$ over the evolution semiflow $\varphi$.
Hence, an evolution operator on $V$ is generating a skew-evolution
semiflow on $Y$.
\end{example}

\begin{example}\rm\label{shift}
If $C=(\varphi, \Phi)$ denotes a skew-evolution semiflow and
$\alpha \in \mathbb{R}$ a parameter, then $C_{\alpha}=(\varphi,
\Phi_{\alpha})$, where
\begin{equation}
\Phi_{\alpha}: T\times  X\rightarrow \mathcal{B}(V), \
\Phi_{\alpha}(t,t_{0},x)=e^{\alpha(t-t_{0})}\Phi(t,t_{0},x),
\end{equation}
is also a skew-evolution semiflow, being the
\emph{$\alpha$-shifted skew-evolution semiflow}.
\end{example}

Other examples of skew-evolution semiflows are given in
\cite{StMe_NA}.

\section{Nonuniform exponential stability}\label{Def_stab}

In this section we define five concepts of exponential stability
for skew-evolution semiflows. For each, an equivalent definition
is given. Also, we will establish some connections between these
concepts and we will emphasize that they are not equivalent.

We will begin by considering the notion of uniform exponential
stability for skew-evolution semiflows, as given in \cite{StMe_NA}
and which was characterized for evolution operators in
\cite{MeSt_IEOT}.

\begin{definition} \label{ues}
A skew-evolution semiflow $C =(\varphi,\Phi)$ is \emph{uniformly
exponentially stable} $(u.e.s.)$ if there exist some constants
$N\geq 1$ and $\alpha>0$ such that, for all $(t,s),(s,t_{0})\in
T$, following relation holds:
\begin{equation}
\left\Vert \Phi(t,t_{0},x)v\right\Vert \leq
Ne^{-(t-s)\alpha}\left\Vert \Phi(s,t_{0},x)v\right\Vert, \ \forall
\ (x,v)\in Y.
\end{equation}%
\end{definition}

An equivalent definition is given by

\begin{remark}\rm
The skew-evolution semiflow $C =(\varphi,\Phi)$ is uniformly
exponentially stable iff there exist $N\geq 1$ and $\alpha>0$ such
that, for all $(t,s)\in  T$, the relation holds:
\begin{equation}
\left\Vert \Phi(t,s,x)v\right\Vert \leq Ne^{-
(t-s)\alpha}\left\Vert v\right\Vert, \ \forall \ (x,v)\in Y.
\end{equation}%
\end{remark}

The nonuniform exponential stability is defined by

\begin{definition} \label{es}
A skew-evolution semiflow $C =(\varphi,\Phi)$ is
\emph{exponentially stable} $(e.s.)$ if there exist a mapping
$N:\mathbb{R}_{+}\rightarrow[1,\infty)$ and a constant $\alpha>0$
such that, for all $(t,s)\in T$, following relation takes place:
\begin{equation}
\left\Vert \Phi(t,t_{0},x)v\right\Vert \leq N(s)e^{-\alpha
t}\left\Vert \Phi(s,t_{0},x)v\right\Vert,\ \forall (x,v)\in Y.
\end{equation}%
\end{definition}

Instead of the previous definition we have

\begin{remark}\rm
The skew-evolution semiflow $C =(\varphi,\Phi)$ is exponentially
stable iff there exist $N\geq 1$ and $\alpha>0$ such that, for all
$(t,s)\in  T$, we have:
\begin{equation}
\left\Vert \Phi(t,s,x)v\right\Vert \leq N(s)e^{-\alpha
t}\left\Vert v\right\Vert,\ \forall \ (x,v)\in Y.
\end{equation}%
\end{remark}

A concept of nonuniform exponential stability for evolution
equations is given by L. Barreira and C. Valls in \cite{BaVa_LNM},
which we will generalize for skew-evolution semiflows. In what
follows, allow us to name this asymptotic propositionerty
\emph{"Barreira-Valls exponential stability".}

\begin{definition}\label{BVes}
A skew-evolution semiflow $C =(\varphi,\Phi)$ is
\emph{Barreira-Valls exponentially stable} $(BV.e.s.)$ if there
exist some constants $N\geq 1$, $\alpha>0$ and $\beta\geq\alpha$
such that, for all $(t,s),(s,t_{0})\in T$, the relation holds:
\begin{equation}
\left\Vert \Phi(t,t_{0},x)v\right\Vert \leq Ne^{-\alpha t}e^{\beta
s}\left\Vert \Phi(s,t_{0},x)v\right\Vert, \ \forall \ (x,v)\in Y.
\end{equation}%
\end{definition}

We also have, as an equivalent definition, the next

\begin{remark}\rm
A skew-evolution semiflow $C =(\varphi,\Phi)$ is Barreira-Valls
exponentially stable iff there some constants $N\geq 1$,
$\alpha>0$ and $\beta\geq\alpha$ such that, for all $(t,s)\in T$,
following relation is verified:
\begin{equation}
\left\Vert \Phi(t,s,x)v\right\Vert \leq Ne^{-\alpha t}e^{\beta
s}\left\Vert v\right\Vert, \ \forall \ (x,v)\in Y.
\end{equation}%
\end{remark}

The asymptotic propositionerty of nonuniform stability is
considered in

\begin{definition} \label{s}
A skew-evolution semiflow $C =(\varphi,\Phi)$ is \emph{stable}
$(s.)$ if there exists a mapping
$N:\mathbb{R}_{+}\rightarrow[1,\infty)$ such that, for all
$(t,s),(s,t_{0})\in  T$, the relation is true:
\begin{equation}
\left\Vert \Phi(t,t_{0},x)v\right\Vert \leq N(s)\left\Vert
\Phi(s,t_{0},x)v\right\Vert, \ \forall \ (x,v)\in Y.
\end{equation}%
\end{definition}

We also have

\begin{remark}\rm
The skew-evolution semiflow $C =(\varphi,\Phi)$ is stable iff
there exists a mapping $N:\mathbb{R}_{+}\rightarrow[1,\infty)$
such that, for all $(t,s)\in  T$, the relation is verified:
\begin{equation}
\left\Vert \Phi(t,s,x)v\right\Vert \leq N(s)\left\Vert
v\right\Vert, \ \forall \ (x,v)\in Y.
\end{equation}%
\end{remark}

Let us remind the propositionerty of exponential growth for
skew-evolution semiflows, given by

\begin{definition} \label{eg}
A skew-evolution semiflow $C =(\varphi,\Phi)$ has
\emph{exponential growth} $(e.g.)$ if there exist two
nondecreasing mappings
$M,\omega:\mathbb{R}_{+}\rightarrow[1,\infty)$ such that, for all
$(t,s),(s,t_{0})\in  T$, we have:
\begin{equation}
\left\Vert \Phi(t,t_{0},x)v\right\Vert \leq
M(s)e^{\omega(t-s)}\left\Vert \Phi(s,t_{0},x)v\right\Vert,\
\forall \ (x,v)\in Y.
\end{equation}%
\end{definition}

Similarly, we have

\begin{remark}\rm
The skew-evolution semiflow $C =(\varphi,\Phi)$ has exponential
growth iff there exist two nondecreasing mappings
$M,\omega:\mathbb{R}_{+}\rightarrow[1,\infty)$ such that, for all
$(t,s)\in  T$, the relation holds:
\begin{equation}
\left\Vert \Phi(t,s,x)v\right\Vert \leq
M(s)e^{\omega(t-s)}\left\Vert v\right\Vert,\ \forall \ (x,v)\in Y.
\end{equation}%
\end{remark}

We obtain following relations concerning the previously defined
asymptotic propositionerties for skew-evolution semiflows.

\begin{remark}\rm From the previous definitions it follows that:
\begin{equation}
(u.e.s.)\Longrightarrow (BV.e.s.)\Longrightarrow
(e.s.)\Longrightarrow(s.)\Longrightarrow (e.g.)
\end{equation}
\end{remark}

The reciprocal statements are not true, as shown in what follows.

\vspace{3mm}

The next example emphasizes a skew-evolution semiflow which is
Barreira-Valls exponentially stable but is not uniformly
exponentially stable.

\begin{example}\rm
Let $ X=\mathbb{R}_{+}$ and $V=\mathbb{R}$. The mapping $\varphi:
T\times\mathbb{R}_{+} \rightarrow \mathbb{R}_{+}$, where $
\varphi(t,s,x)=t-s+x$ \noindent is an evolution semiflow on
$\mathbb{R}_{+}$.

We will consider the function
$u:\mathbb{R}_{+}\rightarrow\mathbb{R}$, given by
$u(t)=e^{2t-t\sin t}$. We define
$$\Phi_{u}(t,s,x)v=\frac{u(s)}{u(t)}v, \ \textrm{with} \ (t,s)\in T,\ (x,v)\in Y.$$
As we have
$$\left| \Phi_{u}(t,s,x)v\right| \leq |v|\cdot e^{t\sin t-s\sin s+2s-2t}\leq |v|e^{3s-2t}=e^{-2t}e^{3t}|v|,$$
for all $(t,s,x,v)\in T\times Y$. It follows that
$C_{u}=(\varphi,\Phi_{u})$ is Barreira-Valls exponentially stable.

Let us suppose now that the skew-evolution semiflow
$C_{u}=(\varphi,\Phi_{u})$ is uniformly exponentially stable.
According to Definition \ref{ues}, there exist $N\geq 1$,
$\alpha>0$ and $t_{1}>0$ such that
$$e^{t\sin t-s\sin s+2s-2t}\leq Ne^{\alpha(s-t)}, \ \forall t\geq s\geq t_{1}.$$
If we consider $t=2n\pi+\frac{\pi}{2}$ and $s=2n\pi$, we have that
$$\exp\left(2n\pi-\frac{3\pi}{2}\right)\leq N\exp\left(-\frac{\pi}{2}\right),$$
which, for $n\rightarrow \infty$, leads to a contradiction, which
proves that $C_{u}$ is not uniformly exponentially stable.
\end{example}

The following example presents a skew-evolution semiflow which is
exponentially stable but not Barreira-Valls exponentially stable.

\begin{example}\rm
Let $ X=\mathbb{R}_{+}$. The mapping $\varphi:
T\times\mathbb{R}_{+} \rightarrow \mathbb{R}_{+}, \
\varphi(t,s,x)=x$ \noindent is an evolution semiflow on
$\mathbb{R}_{+}$.

Let us consider a continuous function
$u:\mathbb{R}_{+}\rightarrow[1,\infty)$ with
$$u(n)=n\cdot 2^{2n}\ \textrm{and} \ u\left(n+\frac{1}{2^{2n}}\right)=1.$$
We define
$$\Phi_{u}(t,s,x)v=\frac{u(s)e^{s}}{u(t)e^{t}}v, \ \textrm{where} \ (t,s)\in T,\ (x,v)\in Y.$$
As following relation
$$\left\Vert \Phi_{u}(t,s,x)v\right\Vert \leq u(s)e^{s}e^{-t}\left\Vert v\right\Vert$$
holds for all $(t,s,x,v)\in T\times Y$, it results that the
skew-evolution semiflow $C_{u}=(\varphi,\Phi_{u})$ is
exponentially stable.

Let us now suppose that the skew-evolution semiflow
$C_{u}=(\varphi,\Phi_{u})$ is Barreira-Valls exponentially stable.
Then, according to Definition \ref{BVes}, there exist $N\geq 1$,
$\alpha>0$, $\beta>0$ and $t_{1}>0$ such that
$$\frac{u(s)e^{s}}{u(t)e^{t}}\leq Ne^{-\alpha t}e^{\beta s}, \ \forall t\geq s\geq t_{1}.$$
For $t=n+\frac{1}{2^{2n}}$ and $s=n$ it follows that
$$e^{n\left(2^{2n}+1\right)}\leq Ne^{n+\frac{1}{2^{2n}}}e^{-\alpha\left(n+\frac{1}{2^{2n}}\right)}e^{\beta n},$$
which is equivalent with $$e^{n\left(2^{2n}-\beta\right)}\leq
Ne^{\frac{1}{2^{2n}}-\alpha\left(n+\frac{1}{2^{2n}}\right)}.$$ For
$n\rightarrow \infty$, a contradiction is obtained, which proves
that $C_{u}$ is not Barreira-Valls exponentially stable.
\end{example}

There exist skew-evolution semiflows that are stable but not
exponentially stable, as results from the following

\begin{example}\rm
Let us consider $ X=\mathbb{R}_{+}$, $V=\mathbb{R}$ and
$$u:\mathbb{R}_{+}\rightarrow[1,\infty)\ \textrm{with the propositionerty}
\ \underset{t\rightarrow\infty}\lim\frac{u(t)}{e^{t}}=0.$$ The
mapping
\[
\Phi_{u}: T\times \mathbb{R}_{+} \rightarrow
\mathcal{B}(\mathbb{R}), \ \Phi_{u}(t,s,x)v=\frac{u(s)}{u(t)}v
\]
is an evolution cocycle. As $|\Phi(t,s,x)v|\leq u(s)|v|$, $\forall
(t,s,x,v)\in T\times Y$, it follows that $C_{u}=(\varphi,
\Phi_{u})$ is a stable skew-evolution semiflow, for every
evolution semiflow $\varphi$ on $\mathbb{R}_{+}$.

On the other hand, if we suppose that $C_{u}$ is exponentially
stable, according to Definition \ref{es}, there exist a mapping
$N:\mathbb{R}_{+}\rightarrow[1,\infty)$ and a constant $\alpha>0$
such that, for all $(t,s),(s,t_{0})\in T$, we have
$$
\left\Vert \Phi(t,t_{0},x)v\right\Vert \leq N(s)e^{-\alpha
t}\left\Vert \Phi(s,t_{0},x)v\right\Vert,\ \forall \ (x,v)\in Y.
$$
It follows that
\[
\frac{u(s)}{N(s)}\leq \frac{u(t)}{e^{\alpha t}}.
\]
For $t\rightarrow\infty$ we obtain a contradiction, and, hence,
$C_{u}$ is not exponentially stable.
\end{example}

Following example gives a skew-evolution semiflow that has
exponential growth but is not stable.

\begin{example}\rm
We consider $ X=\mathbb{R}_{+}$, $V=\mathbb{R}$ and
$$u:\mathbb{R}_{+}\rightarrow[1,\infty)\ \textrm{with the propositionerty}
\ \underset{t\rightarrow\infty}\lim\frac{e^{t}}{u(t)}=\infty.$$
The mapping
\[
\Phi_{u}: T\times \mathbb{R}_{+} \rightarrow
\mathcal{B}(\mathbb{R}), \
\Phi_{u}(t,s,x)v=\frac{u(s)e^{t}}{u(t)e^{s}}v
\]
is an evolution cocycle. We have $|\Phi(t,s,x)v|\leq
u(s)e^{t-s}|v|$, $\forall (t,s,x,v)\in T\times Y$. Hence,
$C_{u}=(\varphi, \Phi_{u})$ is a skew-evolution semiflow, over
every evolution semiflow $\varphi$, and has exponential growth.

Let us suppose that $C_{u}$ is stable. According to Definition
\ref{s}, there exists a mapping
$N:\mathbb{R}_{+}\rightarrow[1,\infty)$ such that $u(s)e^{t}\leq
N(s)u(t)e^{s}$, for all $(t,s)\in T$. If $t\rightarrow\infty$, a
contradiction is obtained. Hence, $C_{u}$ is not stable.
\end{example}

\section{Datko type theorems for the nonuniform exponential
stability}\label{Th_D_R}

A different type of stability for skew-evolution semiflows in the
nonuniform setting is presented in this section, as well a
particular class of skew-evolution semiflows, which allows
connections between various stability types.

\begin{definition}\label{is}
A skew-evolution semiflow $C=(\varphi,\Phi)$ is called
\emph{integrally stable} $(i.s.)$ if there exists a mapping
$D:\mathbb{R}_{+}\rightarrow\mathbb{R}_{+}^{*}$ such that:
\begin{equation}
\int^{\infty}_{s}\left\Vert \Phi(t,t_{0},x)v\right\Vert dt\leq
D(s)\left\Vert \Phi(s,t_{0},x)v\right\Vert,
\end{equation}
for all $(s,t_{0})\in T$ and all $(x,v)\in Y$.
\end{definition}

An equivalent definition can be considered the next

\begin{remark}\rm\label{rem_is}
A skew-evolution semiflow $C=(\varphi,\Phi)$ is integrally stable
iff there exists a mapping
$D:\mathbb{R}_{+}\rightarrow\mathbb{R}_{+}^{*}$ such that:
\begin{equation}
\int^{\infty}_{s}\left\Vert \Phi(t,s,x)v\right\Vert dt\leq
D(s)\left\Vert v\right\Vert,
\end{equation}
for all $s\in \mathbb{R}_{+}$ and all $(x,v)\in Y$.
\end{remark}

\begin{definition}
A skew-evolution semiflow $C=(\varphi,\Phi)$ has \emph{bounded
exponential growth} if $C$ has exponential growth and function $M$
from Definition \ref{eg} is bounded.
\end{definition}

\begin{proposition}\label{proposition_stab}
An integrally stable skew-evolution semiflow $C=(\varphi,\Phi)$
with bounded exponential growth is stable.
\end{proposition}

\begin{proof}
Let us denote $M=\underset{t\geq 0}\sup M(t)$ and
$c=\int_{0}^{1}e^{-\omega(t)}$, where functions $M$ and $\omega$
are given by Definition \ref{eg}.

We observe that for $t\geq s+1$ we have
$$c\leq \int_{0}^{t-s}e^{-\omega(r)dr}=\int_{s}^{t}e^{-\omega(t-\tau)}d\tau$$
and, further,
$$c|<v^{*},\Phi(t,s,x)v>|\leq \int_{s}^{t}e^{-\omega(t-\tau)d\tau}\left\Vert \Phi(t,\tau,\varphi(\tau,s,x))^{*}v^{*}\right\Vert
\left\Vert \Phi(\tau,s,x)v\right\Vert d\tau\leq$$
$$\leq M\int_{s}^{t}\left\Vert \Phi(\tau,s,x)v\right\Vert
d\tau\leq MD(s)\left\Vert v\right\Vert\left\Vert
v\right\Vert^{*},$$ for all $(t,t_{0})\in T$, all $(x,v)\in Y$ and
all $v^{*}\in V^{*}$, function $D$ being given by Remark
\ref{rem_is}.

By taking supremum relative to $\left\Vert v\right\Vert^{*}\leq
1$, we obtain $$\left\Vert \Phi(t,s,x)v\right\Vert\leq
\frac{MD(s)}{c}, \ \forall  t\geq s+1, \ (x,v)\in Y.$$ Finally, it
follows that $$\left\Vert \Phi(t,s,x)v\right\Vert \leq
N(s)\left\Vert v\right\Vert, \ \forall (t,s)\in T, \ \forall
(x,v)\in Y,$$ where we have denoted
$$N(s)=M\left[\frac{D(s)}{c}+e^{\omega(s)}\right],$$ and which proves that $C$ is
stable.
\end{proof}

\begin{definition}\label{eis}
A skew-evolution semiflow $C=(\varphi,\Phi)$ is said to be
\emph{exponentially integrally stable} $(e.i.s.)$ if there exist a
mapping $D:\mathbb{R}_{+}\rightarrow\mathbb{R}_{+}^{*}$ and a
constant $d>0$ such that following relation:
\begin{equation}
\int^{\infty}_{s}e^{(t-s)d}\left\Vert \Phi(t,s,x)v\right\Vert
dt\leq D(s)\left\Vert \Phi(s,t_{0},x)v\right\Vert,
\end{equation}
holds for all $(s,t_{0})\in T$ and all $(x,v)\in Y$.
\end{definition}

We also have

\begin{remark}\rm
A skew-evolution semiflow $C=(\varphi,\Phi)$ is exponentially
integrally stable iff there exist a mapping
$D:\mathbb{R}_{+}\rightarrow\mathbb{R}_{+}^{*}$ and a constant
$d>0$ such that:
\begin{equation}
\int^{\infty}_{s}e^{(t-s)d}\left\Vert \Phi(t,s,x)v\right\Vert
dt\leq D(s)\left\Vert v\right\Vert,
\end{equation}
for all $(t,s)\in T$ and all $(x,v)\in Y$.
\end{remark}

\begin{remark}\rm As a connection between the presented asymptotic propositionerties, we
have:
\begin{equation}
(e.i.s.)\Longrightarrow (i.s.)
\end{equation}
\end{remark}

In what follows, we will emphasize some characterizations of the
various types of nonuniform stability considered in Section
\ref{Def_stab}. We will begin this section by considering a
particular class of skew-evolution semiflows, given in

\begin{definition}
A skew-evolution semiflow $C=(\varphi,\Phi)$ is said to be
\emph{strongly measurable} $(s.m.)$ if for all
$(t_{0},x,v)\in\mathbb{R}_{+}\times Y$ the mapping
$s\mapsto\left\Vert\Phi(s,t_{0},x)v\right\Vert$ is measurable on
$[t_{0},\infty)$.
\end{definition}

\begin{theorem}\label{Datko}
A strongly measurable skew-evolution semiflow $C=(\varphi,\Phi)$
with bounded exponential growth is exponentially stable if and
only if it is exponentially integrally stable.
\end{theorem}

\begin{proof}
\emph{Necessity.} It is a simple verification for
$$d=\frac{\alpha}{2} \ \textrm{and} \ D(t)=\frac{N(t)}{\alpha}, \
t\geq 0.$$ \emph{Sufficiency.} If $C=(\varphi,\Phi)$ is
exponentially integrally stable, then there exists a constant
$d>0$ such that the $d$-shifted skew-evolution semiflow $C_{d}$,
given as in Example \ref{shift}, is integrally stable with bounded
exponential growth.

According to Proposition \ref{proposition_stab}, it follows that
$C_{d}$ is stable, which assures the existence of a mapping
$N:\mathbb{R}_{+}\rightarrow[1,\infty)$ with
$$
\left\Vert \Phi(t,s,x)v\right\Vert \leq N(s)e^{-(t-s)d}\left\Vert
v\right\Vert,\ \forall (t,s)\in T, \ \forall (x,v)\in Y,$$ which
proves that $C$ is exponentially stable.
\end{proof}

\begin{remark}\rm
Theorem \ref{Datko} can be viewed as a Datko type theorem for the
propositionerty of nonuniform exponential stability for
skew-evolution semiflows. The case of uniform stability was
considered in \cite{StMe_NA}. For the particular case of evolution
operators, this result was proved by R. Datko in \cite{Da_JMA} in
the uniform setting and by C. Bu\c{s}e in \cite{Bu_RSMUPT} for the
nonuniform case.
\end{remark}

Let us denote by $\mathcal{F}$ the set of all nondecreasing
functions $F:\mathbb{R}_{+}\rightarrow\mathbb{R}_{+}$ with the
propositionerties $F(0)=0$ and $F(t)>0$, $\forall t>0$.

\begin{remark}\rm
Analogously to the uniform case studied in \cite{StMe_NA}, the
proof of Theorem \ref{Datko} can be easily adapted to prove a
variant of Rolewicz type for the propositionerty of exponential
stability of skew-evolution semiflows in the nonuniform setting,
as given by
\end{remark}

\begin{theorem}\label{Rolewicz}
Let $C=(\varphi,\Phi)$ be a strongly measurable skew-evolution
semiflow with exponential growth. Then $C$ is exponentially stable
if and only if there exist two mappings
$F,R:\mathbb{R}_{+}\rightarrow\mathbb{R}_{+}$ and a constant $d>0$
with $F\in\mathcal{F}$ and:
\begin{equation}
\int^{\infty}_{s}F\left(e^{(t-s)d}\left\Vert
\Phi(t,s,x)v\right\Vert dt\right)\leq R(s)F\left(\left\Vert
v\right\Vert\right),
\end{equation}
for all $(s,x,v)\in\mathbb{R}_{+}\times Y$.
\end{theorem}

\begin{remark}\rm
For the particular case of evolution operators, Theorem
\ref{Rolewicz} was proved by S. Rolewicz in \cite{Ro_JMAA} for the
propositionerty of uniform exponential stability.
\end{remark}

\begin{remark}\rm
By means of the methods used in the proofs of Proposition
\ref{proposition_stab} and of Theorem \ref{Datko}, one can obtain
a Datko type theorem for the exponential stability of
Barreira-Valls type, in the case of skew-evolution semiflows in
the nonuniform setting, as shown by
\end{remark}

\begin{theorem}\label{Datko_BV}
Let $C=(\varphi,\Phi)$ be a strongly measurable skew-evolution
semiflow with exponential growth. Then $C$ is Barreira-Valls
exponentially stable if and only if there exist some constants
$N\geq 1$, $a>0$ and $b\geq a$ such that:
\begin{equation}
\int^{\infty}_{s}e^{at}\left\Vert \Phi(t,s,x)v\right\Vert dt\leq
Ne^{bs}\left\Vert v\right\Vert,
\end{equation}
for all $(s,x,v)\in\mathbb{R}_{+}\times Y$.
\end{theorem}

\begin{remark}\rm
Analogously, a Rolewicz type theorem can be given for the
propositionerty of Barreira-Valls exponential stability, in the
case of skew-evolution semiflows.
\end{remark}

\section{A Barbashin type theorem for the nonuniform exponential
stability}\label{Th_B}

In this section let us consider a particular class of
skew-evolution semiflows, given by

\begin{definition}
A skew-evolution semiflow $C=(\varphi,\Phi)$ is said to be
\emph{$*$-strongly measurable} $(\ast-s.m.)$ if for every
$(t,t_{0},x,v^{*})\in T\times X\times V^{*}$ the mapping defined
by
$s\mapsto\left\Vert\Phi(t,s,\varphi(s,t_{0},x))^{*}v^{*}\right\Vert$
is measurable on $[t_{0},t]$.
\end{definition}

The main result of this section is

\begin{theorem}\label{Barbashin}
Let $C=(\varphi,\Phi)$ be a $*$-strongly measurable skew-evolution
semiflow with exponential growth. If there exist a constant $b>0$
and a mapping $B:\mathbb{R}_{+}\rightarrow [1,\infty)$ such that:
\begin{equation}
\int^{t}_{s}e^{(t-\tau)b}\left\Vert
\Phi(t,\tau,\varphi(\tau,s,x))^{*}v^{*}\right\Vert d\tau\leq
B(t)\left\Vert v^{*}\right\Vert,
\end{equation}
for all $(t,s)\in T$ and all $(x,v^{*})\in X\times V^{*}$, then
$C$ is exponentially stable.
\end{theorem}

\begin{proof}
For $t\geq s\geq 0$ we will denote
$$f_{s}(t)=M(s)B(t)e^{tb}e^{\omega(t)} \ \textrm{and} \ K(s)=\int_{0}^{1}\frac{du}{f_{s}(u)},$$ where the functions $M$ and
$\omega$ are given by Definition \ref{eg}.

We remark that, if $t\geq s+1$, then $$K(s)\leq
\int^{t-s}_{0}\frac{du}{f_{s}(u)}=\int^{t}_{s}\frac{d\tau}{f_{s}(\tau-s)}.$$
It follows that $$B(t)e^{(t-s)b}K(s)|<v^{*},\Phi(t,s,x)v>|\leq$$
$$\leq
\int^{t}_{s}\frac{e^{(t-s)b}|<\Phi(t,\tau,\varphi(\tau,s,x))^{*}v^{*},
\Phi(\tau,s,x)v>|}{M(s)e^{(\tau-s)b}e^{\omega(\tau-s)}}d\tau\leq$$
$$\leq \int^{t}_{s}e^{(t-\tau)b}\left\Vert \Phi(t,\tau,\varphi(\tau,s,x))^{*}v^{*}\right\Vert
\left\Vert v\right\Vert d\tau \leq B(t)\left\Vert
v\right\Vert\left\Vert v^{*}\right\Vert,$$ which implies
$$\left\Vert \Phi(t,s,x)v\right\Vert\leq
\frac{e^{-(t-s)b}}{K(s)}\left\Vert v\right\Vert$$ for all $t \geq
s+1$ and all $(x,v)\in Y$.

Now, if we consider $t\in [s,s+1)$, we have $$\left\Vert
\Phi(t,s,x)v\right\Vert\leq M(s)e^{\omega(t-s)}\left\Vert
v\right\Vert\leq M(s)e^{\omega(1)}\left\Vert v\right\Vert\leq
M(s)e^{b+\omega(1)}e^{-b(t-s)}\left\Vert v\right\Vert.$$ Finally,
we obtain, $$\left\Vert \Phi(t,s,x)v\right\Vert\leq
N(s)e^{-(t-s)b}\left\Vert v\right\Vert,$$ for all $(t,s)\in T$ and
all $(x,v)\in X\times V$, where we have denoted
$$N(s)=M(s)e^{b+\omega(1)}+\frac{1}{K(s)},$$ and which proves the
exponential stability of the skew-evolution semiflow $C$.
\end{proof}

\begin{remark}\rm
Theorem \ref{Barbashin} is a generalization of a known result of
E.A. Barbashin emphasized in \cite{Ba_NAU}. A similar result was
obtained by C. Bu\c{s}e, M. Megan, M. Prajea and P. Preda for the
uniform exponential stability of evolution operators in
\cite{BuMePrPr_IEOT}.
\end{remark}

\begin{remark}\rm
Analogously as in the proof of Theorem \ref{Barbashin}, one can
prove a Barbashin type theorem for the propositionerty of
Barreira-Valls exponential stability, in the case of
skew-evolution semiflows.
\end{remark}

\textbf{Acknowledgments.} This work is financially supported from
the Exploratory Research Grant CNCSIS PN II ID 1080 No. 508/2009
of the Romanian Ministry of Education, Research and Innovation.

{\footnotesize

\vspace{5mm}

\noindent\begin{tabular}[t]{ll}

(Codru\c{t}a Stoica) \textsc{Department of Mathematics and
Computer Science}, \\ \textsc{"Aurel Vlaicu"} \textsc{University
of Arad}, \textsc{Romania}\\
\textit{E-mail address:}
\texttt{codruta.stoica@math.u-bordeaux1.fr}

\vspace{3mm}\\

(Mihail Megan) \textsc{Academy of Romanian Scientists, Bucharest,
Romania}\\ \textsc{Faculty of Mathematics and Computer Science}\\
\textsc{West University of Timi\c{s}oara}, \textsc{Romania}\\
\textit{E-mail address:} \texttt{megan@math.uvt.ro}
\end{tabular}
}

\end{document}